\documentstyle[12pt]{article}
\title{ The relation between the
decomposition  of comodules and coalgebras
 \thanks {This work was supported by the  National Natural Science Foundation  }}
\author{
Shouchuan Zhang  \\ Department  of Mathematics, Hunan University\\
Changsha  410082, \
 P.R.China. \ \
E-mail:z9491@yahoo.com.cn\\
}
\date{}
\begin{document}
% \tableofcontents
\newtheorem{Theorem}{\quad Theorem}[section]
\newtheorem{Proposition}[Theorem]{\quad Proposition}
\newtheorem{Definition}[Theorem]{\quad Definition}
\newtheorem{Corollary}[Theorem]{\quad Corollary}
\newtheorem{Lemma}[Theorem]{\quad Lemma}
\newtheorem{Example}[Theorem]{\quad Example}
\maketitle \addtocounter{section}{-1}

 \begin {abstract}   T. Shudo and H. Miyamito \cite{SM78} showed that
  $C$ can  be decomposed
into a direct sum  of its indecomposable subcoalgebras of $C$.
 Y.H. Xu \cite {XF92}  showed that the decomposition  was unique.
He also  showed that
  $M$ can uniquely  be decomposed
into a direct sum  of the weak-closed indecomposable subcomodules
of $M$(we call the decomposition
 the weak-closed indecomposable
decomposition ) in  \cite{XSF94}. In this paper, we   give
 the relation between the
   two decomposition.
   We show that if $M$  is a full, $W$-relational
   hereditary $C$-comodule, then the following conclusions hold:

(1) $M$ is indecomposable iff $C$ is indecomposable;

(2) $M$ is relative-irreducible iff $C$ is irreducible;

(3) $M$ can be decomposed into a direct sum of
 its weak-closed  relative-irreducible subcomodules iff
$C$ can be decomposed into a direct sum of its irreducible
subcoalgebras.     \\
We also obtain the relation between  coradical  of $C$- comodule
$M$ and radical of algebra  $C(M)^*$

 \end {abstract}

\section {Introduction and Preliminaries}

The decomposition of coalgebras and comodules is an important
subject in study of Hopf algebras.
  T. Shudo and H. Miyamito \cite{SM78} showed that
  $C$ can  be decomposed
into a direct sum  of its indecomposable subcoalgebras of $C$.
 Y.H. Xu \cite {XF92}  showed that the decomposition  was unique.
He also  showed that
  $M$ can uniquely  be decomposed
into a direct sum  of the weak-closed indecomposable subcomodules
of $M$(we call the decomposition
 the weak-closed indecomposable
decomposition ) in  \cite{XSF94}. In this paper, we   give
 the relation between the
   two decomposition.
   We show that if $M$  is a full, $W$-relational
   hereditary $C$-comodule, then the following conclusions hold:

(1) $M$ is indecomposable iff $C$ is indecomposable;

(2) $M$ is relative-irreducible iff $C$ is irreducible;

(3) $M$ can be decomposed into a direct sum of
 its weak-closed  relative-irreducible subcomodules iff
$C$ can be decomposed into a direct sum of its irreducible
subcoalgebras.     \\
We also obtain the relation between  coradical  of $C$- comodule
$M$ and radical of algebra  $C(M)^*$

Let $k$ be a field, $M$ be a $C$-comodule, $N$ be a subcomodule of
$M$, $E$ be a subcoalgebras of  $C$ and $P$ be an ideal of  $C^*$.
As in  \cite{XSF94},  we define:

$E^{ \bot C^* } = \{ f \in C^* \mid  f(E) = 0 \}$.

$P^{ \bot C } = \{ c \in C \mid  P(c) = 0 \}$.

$N^{ \bot C^* } = \{ f \in C^* \mid  f \cdot N = 0 \}$.

$P^{ \bot M } = \{ x \in M \mid  P \cdot x = 0 \}$.

Let $<N>$ denote  $N^{ \bot C^*  \bot M } $.
 $<N>$ is called the closure of  $N$. If  $<N> = N$,
then  $N$ is called closed. If  $N  = C^*x$, then $<N>$ is denoted
by $<x>$. If  $<x> \subseteq N $ for any $x \in N$, then $N$ is
called weak-closed. It is clear that any closed subcomodule is
weak-closed. If  $\rho (N) \subseteq N \otimes E  $, then $N$ is
called an $E$-subcomodule of $M$. Let
$$M_E = \sum \{ N \mid N \hbox {  is  a  subcomodule  of } M
\hbox { and } \rho (N) \subseteq N \otimes E \}. $$ We call  $M_E$
a component of $M$ over $E$. If $M_E$ is some component of $M$ and
$M_F \subseteq M_E$ always implies  $M_F = M_E$ for any non-zero
component $M_F$,  then  $M_E$ is called a minimal component of
$M$.

Let $\{m_{ \lambda }  \mid  \lambda \in \Lambda \}$ be the basis
of  $M$ and $C(M)$ denote the subspace of $C$ spanned by

$W(M) = \{ c \in C \mid  {~there~ exists~ an~} m \in M
 {~with~}
\rho(m) = \sum m_{ \lambda } \otimes c_{\lambda } {~such ~that~ }
 c_{\lambda _0 } = c  {~for ~some~}   \lambda _0 \in \Lambda  \}$.\\
E. Abe  in \cite[P129]{Ab80} checked that  $C(M)$ is a
subcoalgebra of  $C$. It is easy to know that if $E$ is
subcoalgebra of $C$ and $ \rho (M) \subseteq M \otimes E$, then
$C(M) \subseteq E$. If $C(M) = C$, then $M$ is called a full
$C$-comodule. If $D$ is a simple subcoalgbra  of   $C$ and  $M_D
\not= 0$ or  $D=0$, then $D$ is called faithful to  $M$. If  every
simple subcoalgbra of $C$ is faithful to  $M$, then  $M$ is called
a component faithfulness $C$-comodule.

     Let     $X$ and $Y$  be subspaces of coalgebra $C$.
     Define $X \wedge Y$  to be the kernel of the composite
     $$C \stackrel {\Delta } {\rightarrow } C \otimes C \rightarrow
     C/X \otimes C/Y. $$
     $X \wedge Y$  is called a wedge of $X$  and $Y$.

\section{The relation between the
decomposition  of comodules and coalgebras}
\begin{Lemma} \label{7.1.1}
 Let $N$ and $L$ be subcomodules of  $M$.
Let $D$ and $E$ be subcoalgebras of  $C$. Then

(1)  $N^{ \bot C^* } = (C(N))^{ \bot C^* }$;

(2)  $ M_E = M_{C(M_E)}$;

(3) $N$ is closed iff there exists a subcoalgebra  $E$ such that
$N =  M_E $;

(4)  $M_D \not= 0$ iff $D^{\bot C^* \bot M} \not= 0 $;

(5)  If $D \cap C(M) = 0$, then $M_D=0$;

(6)  If $D$ is a simple subcoalgebra of  $C$, then
$C(M_D) =   \left \{  \begin{array}{l l} D & { \rm if~} M_D \not= 0\\
0 & {\rm if~ }  M_D = 0
\end{array} \right.$;

(7) If $D\cap E=0$, then $M_D\cap M_E = 0$;

(8)  If  $D$ and $E$ are simple subcoalgebras and  $M_D \cap M_E =
0$ with  $M_D \not= 0$ or  $M_E \not= 0$,
 then  $D \cap E = 0$.

(9)  If $M_E$ is the minimal component of $M$ and  $0 \not= N
\subseteq M_E$,
 then  $C(M_E) = C(N)$.

\end{Lemma}
{\bf Proof}.
 (1)
Let $\{m_{\lambda} \mid \lambda\in\Lambda\}$ be a basis of $N$.
 For any  $f\in N^{\bot C^*}$, if $c\in W(N)$, then there exists
$m\in N$ with $\rho(m)= \sum_{\lambda\in\Lambda}m_{\lambda}\otimes
c_{\lambda}$ such that
 $c_{\lambda_0}=c$ for some $\lambda_0 \in\Lambda$.
  Since $f\cdot m=\sum m_{\lambda}
f(c_{\lambda})=0$, $f(c_{\lambda})=0$ for any $\lambda\in\Lambda$.
 Obviously, $f(c)=f(c_{\lambda _0})=0$.
Considering  that $C(N)$ is the space spanned by $W(N)$, we have
$f\in C(N)^{\bot C^*}$. Conversely, if $f\in C(N)^{\bot C^*}$,
then $f\cdot m=\sum m_{\lambda}f(c_{\lambda})=0$ for any $m\in N$,
i.e. $f\in N^{\bot C^*}$. This shows that $N^{\bot C^*}=C(N)^{\bot
C^*}$.

(2) Since $M_E$ is an $E$-subcomodule of $M$, $C(M_E)\subseteq E$
and $M_{ C(M_E)}\subseteq M_E$. Since $M_E$ is a
$C(M_E)$-subcomodule, $M_E \subseteq M_{C(M_E)}$.

(3) If $N$ is closed, then $N^{\bot C^*\bot M}=N$. Let $E=C(N)$.
Obviously, $N\subseteq M_E$. By Lemma \ref{7.1.1}(1),
 $N^{ \bot C^{\star}}\cdot M_E=E^{\bot C^{\star}}
\cdot M_E=0$. Thus
 $M_E\subseteq N^{\bot C^{\star}\bot M}$ and $M_E\subseteq N$.
This shows that $M_E=N$. Conversely, if $M_E=N$,  obviously,
$(M_E) \subseteq (M_E)^{ \bot C^*  \bot M }.$ Thus it is
sufficient to show that
    $$<M_E> \subseteq  M_E.$$
    Let $L= (M_E)^{\bot C^* \bot M} = <M_E>.$ We see that
  $C(L)^{ \bot C^* } = L^{ \bot C^* }
= (M_E)^{ \bot C^* \bot M \bot C^* } = (M_E)^{ \bot C^* }=
C(M_E)^{ \bot C^* }$. Thus $C(L) = C(M_E)$ and  $ L \subseteq
M_{C(L)} = M _{C(M_E)} = M_E $.

(4)  If  $D^{ \bot C^* \bot M } \not= 0$, let $\{m_{\lambda} \mid
\lambda \in\Lambda \}$ be a basis of $D^{\bot C^* \bot M}$ and
$\rho (x) = \sum m_ \lambda \otimes d_ \lambda $ for any $x \in
D^{\bot C^* \bot M } $.
 Since   $D^{ \bot C^* } \cdot x  =
\sum_{\lambda } m_\lambda  D^{ \bot C^* } (d_ \lambda) = 0$,
$d_\lambda
 \in
D^{ \bot C^* \bot C } = D$ for any $\lambda \in \Lambda$, which
implies
 that
$D^{ \bot C^* \bot M}$  is a $D$-subcomodule. Therefore $ 0 \not=
D^{ \bot C^* \bot M} \subseteq M_D$, i.e. $M_D \not= 0$.
 Conversely, if  $M_D\not= 0$,  then we have that
 \begin {eqnarray*}
 0 \not= M_D &\subseteq & (M_D)^{ \bot C^* \bot M} \\
 &  = & C(M_D)^{ \bot C^* \bot M} \hbox { \ \ (by part (1)) }  \\
& \subseteq & D^{ \bot C^* \bot M}.
\end {eqnarray*}

(5)  Since $M_D \subseteq M$, $C(M_D) \subseteq C(M)$. Obviously,
$C(M_D) \subseteq D$. Thus $C(M_D) \subseteq C(M) \cap D = 0$,
which implies  that   $M_D = 0$.

(6)  If  $M_D = 0$, then $C(M_D) = 0$. If $M_D \not= 0$,
 then $0 \not= C(M_D) \subseteq D$. Since   $D$ is a simple
subcoalgebra, $C(M_D) = D$.

(7)  If $x \in M_D \cap M_E$, then
 $\rho (x) \in (M_D \otimes D) \cap (M_E \otimes E)$,
and
 \begin {eqnarray} \rho (x) =   \sum_{ i = 1 }^{n}  x_i \otimes d_i =
 \sum_{ j = 1 }^{m}  y_j \otimes e_j,
 \label {e.1.1}
 \end {eqnarray}
where $x_i, \cdots , x_n$ is linearly independent and   $d_i \in
D$ and $ e_j \in E$. Let $ f_l \in M^*$ with  $ f_l(x_i) = \delta
_{li}$ for  $ i, l = 1, 2, \cdots n$ . Let  $f_l \otimes id$ act
on equation ({\ref{e.1.1}).
 We have that
$d_l = \sum_{ j = 1 }^{m} f_l(y_j)e_j \in E,$  which implies
  $d_l \in D \cap E = 0$ and  $d_l = 0$  for $ l = 1, \cdots , n$.
 Therefore  $M_D \cap M_E = 0$.

(8)  If $D \cap E \not= 0$, then $D =E$. Thus $M_D = M_E$ and
$M_D \cap M_E = M_E = M_D = 0$.
 We get a contradiction. Therefore  $D \cap E = 0$.

(9) Obviously, $C(N) \subseteq C(M_E)$. Conversely,since  $0 \not=
N
 \subseteq M_{C(N)} \subseteq M_{C(M_E)} = M_E$  and
$M_E$ is a minimal component, $M_{C(N)} = M_E$. By the definition
of component, $C(M_E) \subseteq C(N)$. Thus $C(N) = C(M_E)$.
 $\Box$

\begin{Proposition} \label{7.1.2}
 If $E$ is a subcoalgebra of $C$,then
the following conditions are equivalent.

(1)  $M_E$ is a minimal component of $M$.

(2)  $C(M_E)$ is a simple subcoalgebra of $C$.

(3)  $M_E$ is a minimal closed subcomodule of $M$.
\end{Proposition}
{\bf Proof}.  It is easy to check that $M_E = 0$ iff  $C(M_E)= 0$.
Thus  (1), (2) and (3) are equivalent when $C(M_E) =0$,. We now
assume that $M_E \not= 0$.

(1)  $\Longrightarrow$ (2) Since   $ 0 \not= M_E$, there exists a
non-zero finite dimensional simple subcomodule $N$  of $M$ such
that
 $N \subseteq M_E$. By  \cite[Lemma 1.1]{XSF94},
$N$ is  a simple $C^*$-submodule of $M$. Since $M_E$ is a minimal
component, $C(N) = C(M_E)$  by Lemma \ref {7.1.1} (9). Let $D =
C(N) = C(M_E)$. By Lemma \ref{7.1.1}(1), $(0:N)_{C^*} = N^{ \bot
C^*} = C(N)^{ \bot C^*} = D^{ \bot C^*}$. Thus $N$ is a faithful
simple  $C^*/D^{ \bot C^*}$-module, and so  $C^*/D^{ \bot C^*}$ is
a simple algebra.It is clear that $D$ is a simple subcoalgebra of
$C.$

$(2 ) \Longrightarrow (3)$  If $0\not= N \subseteq M_E$ and  $N$
is a closed subcomodule of  $M$, then by Lemma \ref {7.1.1}(3)
 there exists a subcoalgebra $F$ of  $C$ such that  $N = M_F$.
Since $0 \not= C(N) = C(M_F) \subseteq C(M_E)$ and   $C(M_E)$ is
simple, $C(M_F) = C(M_E)$. By Lemma \ref{7.1.1}(2), $ N = M_F =
M_{C(M_F)} = M_{C(M_E)} = M_E$, which implies  that $M_E$ is a
minimal closed subcomodule.

$(3) \Longrightarrow (1)$, If  $0\not= M_F \subseteq M_E$,
  then $M_F$ is a closed subcomodule by Lemma \ref{7.1.1}(3) and
   $M_E = M_F$, i.e.
$M_E$ is a minimal component.

This completes the proof.
 $\Box$

Let

${ \bf \cal C}_0 =  \{ D  \mid  D$ is a simple subcoalgebra of
$C$ \};

${ \bf \cal C}(M)_0 =  \{ D  \mid  D$ is a simple subcoalgebra of
$C$ and  $D \subseteq C(M)$ \};

${ \bf \cal M}_0 =  \{ N  \mid  N$ is a minimal closed subcomodule
of  $M$ \};

${ \bf \cal C}(M)_1 =  \{ D  \mid  D$ is a faithful simple
subcoalgebra of  $C$ to $M$ \};

$C_0 = \sum  \{ D  \mid  D$ is  a simple subcoalgebra of $C$ \};

$M_0 = \sum \{ N  \mid  N$ is a minimal closed subcomodule of  $M$
\};

$C_0$ and  $M_0$  are called the coradical of coalgebra  $C$ and
the coradical of comodule  $M$ respectively. If $M_0 =M$, then M
is called cosemisimple.  By  Lemma \ref {7.1.1}(5),

${ \bf \cal C}(M)_1 \subseteq { \bf \cal C}(M)_0$

\begin {Theorem}  \label {7.1.3}
$\psi  \left \{  \begin{array}{l} { \bf \cal C}(M)_1  \longrightarrow { \bf \cal M}_0\\
 D \longmapsto M_D
\end{array} \right.$
is bijective.

\end{Theorem}
{\bf Proof}.  By Proposition \ref{7.1.2}, $\psi$ is a map. Let $D$
and $E \in$ { \bf \cal C}$(M)_1$ and $\psi (D) = \psi (E)$, i.e.
$M_D = M_E $. By Lemma \ref {7.1.1} (6), we have that $D =C(D_D)=
C(M_E)=E.$ If $N \in { \bf \cal M}_0$, then $N = M_{C(N)}$ and
$C(N)$ is a simple subcoalgebra by Lemma \ref{7.1.1}(3) and
Proposition \ref{7.1.2}. Thus $\psi (C(N)) = N,$  which implies
that  $ \psi$ is surjective. $\Box$

 In \cite{XSF94} and \cite{XF92}, Xu defined the equivalence relation
for coalgebra and for comodule as follows:

\begin{Definition} \label{7.1.4}
   We say that
 $D \sim E$ for $D$ and  $E \in { \bf \cal C}_0$ iff for any pair of
subclasses ${ \bf \cal C}_1$  and   ${ \bf \cal C}_2$
 of  ${ \bf \cal C}_0$ with
 $D \in    { \bf \cal C}_1$ and $ E \in { \bf \cal C}_2 $ such that
${ \bf \cal C}_1 \cup { \bf \cal C}_2  = { \bf \cal C}_0$ and ${
\bf \cal C}_1 \cap { \bf \cal C}_2 = \emptyset $, there exist
elements $D_1 \in { \bf \cal C}_1$ and $ E_1 \in { \bf \cal C}_2 $
such that
 $D_1 \wedge E_1 \not= E_1 \wedge D_1 $. Let $[D]$ denote the
equivalence class which contains  $D$.

   We say that
 $N \sim L$ for $N$ and  $L \in { \bf \cal M}_0$ iff for any pair of
subclasses ${ \bf \cal M}_1$  and   ${ \bf \cal M}_2$
 of  ${ \bf \cal M}_0$ with
 $N \in    { \bf \cal M}_1$ and $ L \in { \bf \cal M}_2 $ such that
${ \bf \cal M}_1 \cup { \bf \cal M}_2  = { \bf \cal M}_0$ and ${
\bf \cal M}_1 \cap { \bf \cal M}_2 = \emptyset $, there exist
elements $N_1 \in { \bf \cal M}_1$              and $ L_1 \in {
\bf \cal M}_2 $ such that
 $N_1 \wedge L_1 \not= L_1 \wedge N_1 $. Let $[N]$ denote the
equivalence class which contains  $N$.
\end {Definition}

\begin{Definition} \label{7.1.5}

 If $D \wedge E = E \wedge D$ for  any simple subcoalgebras $D$ and $E$
 of $C$,
 then $C$ is called $\pi$-commutative.
  If $N \wedge L = L \wedge N$ for any minimal closed subcomodules
$N$ and $L$ of $M$, then $M$ is called  $\pi$-commutative.
\end {Definition}

Obviously, every cocommutative coalgebra is $\pi$-commutative.
  By \cite[Theorem 3.8 and Theorem 4.18]{XSF94},
$M$ is  $\pi$-commutative iff    $M$ can   be decomposed into a
direct sum  of the weak-closed relative-irreducible
 subcomodules of $M$ iff every equivalence class of $M$ contains only one element.
  By \cite{XF92}, $C$ is  $\pi$-commutative iff
   $C$ can   be decomposed
into a direct sum  of irreducible
 subcoalgebras  of $C$ iff equivalence every class of $C$ contains only one
 element.

\begin{Lemma} \label{7.1.6}          Let
 $D$, $E$ and $F$ be subcoalgebras of  $C$. $N$, $L$ and $T$
 be subcomodules of  $M$. Then

(1)  $M_D \wedge M_F = M_{C(M_D) \wedge C(M_F)} \subseteq M_{D
\wedge F}$;

(2)  If  $D$ and  $E$ are faithful simple subcoalgebras of  $C$
 to $M$, then
$ M_D \wedge M_E =  M_{D \wedge F}$;

(3)   $ M_{D + E} \supseteq  M_D + M_E$;

(4)  If $F = \sum $\{ $D_{\alpha} \mid  \alpha \in \Omega $\} and
$\{ D_{\alpha} \mid  \alpha \in \Omega \}  \subseteq  { \bf \cal
C}_0$, then $M_F = \sum \{ M_{D_{\alpha}} \mid \alpha \in \Omega
\}$.In particular, $M_{C_0} = M_0$.

(5) $  (N + L) \wedge T \supseteq N \wedge T + L \wedge T$;

(6)  $(D + E) \wedge F \supseteq D \wedge F + E \wedge F$.
\end{Lemma}
{\bf Proof}. (1)     We see that
\begin {eqnarray*}
 M_D \wedge M_E &=&
\rho ^{-1}(M \otimes (M_D)^{ \bot C^* \bot C } \wedge
 (M_E)^{ \bot C^* \bot C }) \hbox { \ ( by \cite[Proposition 2.2(1)]{XSF94}) }\\
&=&  \rho ^{-1}(M \otimes C(M_D) \wedge C(M_E)) \hbox { \ (by
Lemma \ref{7.1.1}(1)).}
\end {eqnarray*}
By the definition of component, subcomodule $M_D \wedge M_E
\subseteq M_{C(M_D) \wedge C(M_E)}$. It follows  from the equation
above  that
 $M_{C(M_D) \wedge C(M_E) } \subseteq M_D \wedge M_E$.
Thus
$$M_D \wedge M_E = M_{C(M_D) \wedge C(M_E)}$$
and
 $$M_{C(M_D) \wedge C(M_E) } \subseteq M_{D \wedge E}.$$

(2)  Since $D$ and $E$ are faithful simple subcoalgebras of  $C$
to
 $M$, $C(M_D) = D$ and $C(M_E) = E$  by  Lemma \ref {7.1.1}(6). By  Lemma \ref{7.1.6}(1),
 $M_D \wedge M_E =  M_{D \wedge E}$.

(3) It is trivial.

(4) Obviously $ M_F \supseteq \sum \{ M_{D_{ \alpha }} \mid \alpha
\in \Omega \} $. Conversely, let $N = M_F$. Obviously $N$ is  an
$F$-comodule. For any $x \in N$, let $ L = C^*x$. it is clear
 that  $L$ is a finite dimensional comodule over $F$.
By \cite[Lemma 14.0.1] {Sw69a}, $L$ is a completely reducible
module over
 $F^*$. Thus  $L$ can be decomposed into a direct sum of simple
$F^*$-submodules:
$$L = L_1 \oplus L_2 \oplus  \cdots  \oplus L_n,$$
where  $L_i$ is a simple $F^*$-submodule. By \cite [Proposition
1.16] {XSF94},
 $<L_i> =
(L_i)^{ \bot F^* \bot N}$ is a minimal closed  $F$-subcomodule of
$N$.
 By Theorem \ref{7.1.3}, there exists a simple subcoalgebra
$D_{\alpha _i}$ of $F$ such that
 $<L_i> =
N_{ D_{\alpha_i}}$. Obviously, $N_{ D_{\alpha_i}} \subseteq M_{
D_{\alpha _i}}$. Thus $L \subseteq \sum_{ i = 1}^{n} <L_i> =
\sum_{i = 1}^{ n} N_{D_{\alpha _i}} \subseteq \sum \{ M_{D_
{\alpha}}  \mid \alpha  \in \Omega \}$. Therefore
 $M_F =  \sum \{ M_{ D_{\alpha}}  \mid
\alpha  \in \Omega \}$. If $C_0 = F = \sum \{ D  \mid  D \in { \bf
\cal C}_0 \}$, then $M_{C_0} = M_F = \sum \{ M_D  \mid  D  \in {
\bf \cal C}_0  \}=M_0$ by Theorem \ref{7.1.3},

(5) and (6)  are trivial.$\Box$

\begin{Lemma} \label{7.1.7} Let $N$ be a subcomodule
of  $M$, and let $D$, $E$ and $F$ be simple subcoalgebras of $C$.
Then

(1) $D \sim 0$ iff $D = 0$;  {~}{~} $M_D \sim 0$ iff $M_D = 0$;

We called the equivalence class which contains zero a zero
equivalence class.

(2)  If  $D$ and $E$ are faithful to  $M$ and
  $M_D \sim M_E$, then  $D \sim E$;

(3)  $[M_D] \subseteq \{ M_E \mid  E \in [D] \}$;

(4)  If  $D$ is faithful to  $M$, then
   $[M_D] \subseteq \{ M_E \mid  E \in [D]$ and  $E$  is faithful to  $M$  \}.
   \end{Lemma}
 {\bf Proof}.
(1) If $ D \sim 0$ and $ D \not=  0$, let  ${ \bf \cal C}_1 = \{ 0
\}$ and ${ \bf \cal C}_2 = $\{ $F \mid F \not= 0$, $F \in { \bf
\cal C}_0 $\}. Thus ${ \bf \cal C}_1 \cup { \bf \cal C}_2 =  { \bf
\cal C}_0$ and
 ${ \bf \cal C}_1 \cap { \bf \cal C}_2 = \emptyset $ with
$0 \in  { \bf \cal C}_1$ and $D \in { \bf \cal C}_2$. But for any
$D_1 \in { \bf \cal C}_1$ and $E_1 \in { \bf \cal C}_2$, since
$D_1 = 0$,  $D_1 \wedge E_1 = E_1 \wedge D_1 = E_1$. We get a
contradiction. Thus  $D = 0$. Conversely, if  $D = 0$, obviously $
D \sim 0$. Similarly, we can show that   $M_D \sim 0$ iff $M_D =
0$.

(2) For any pair of subclasses
 ${ \bf \cal C}_1$ and
 ${ \bf \cal C}_2$ of ${\bf \cal C}_0$, if
${ \bf \cal C}_1 \cap { \bf \cal C}_2 = \emptyset $ and   ${ \bf
\cal C}_1  \cup { \bf \cal C}_2 = { \bf \cal C}_0$ with
 $ D \in { \bf \cal C}_1$ and
$E \in { \bf \cal C}_2$, let
 ${ \bf \cal M}_1  = $\{ $M_F \mid F \in { \bf \cal C}_1$  and  $F$
 is faithful to  $M$  \} and
 ${ \bf \cal M}_2  =$ \{ $M_F \mid F \in { \bf \cal C}_2 $ and  $F$
 is faithful to  $M$ \}.
By Theorem \ref{7.1.3},
 ${ \bf \cal M}_0 = { \bf \cal M}_1 \cup { \bf \cal M}_2$ and
 ${ \bf \cal M}_1 \cap { \bf \cal M}_2 =
\emptyset $ . Obviously, $M_D \in  { \bf \cal M}_1$ and $M_E \in {
\bf \cal M}_2$. Since  $ M_D \sim M_E$,there exist  $M_{D_1} \in {
\bf \cal M}_1$ and
 $M_{E_1} \in { \bf \cal M}_2$ such that
 $M_{D_1} \wedge  M_{E_1} \not=   M_{E_1} \wedge  M_{D_1}$, where
$D_1 \in { \cal C}_1$ and  $E_1 \in { \cal C}_2$. By Lemma
\ref{7.1.6}(2), $M_{ D_1 \wedge E_1} \not=  M_{E_1 \wedge D_1}$.
Thus $D_1 \wedge E_1 \not= E_1 \wedge D_1$. Obviously
 $D {~and~}  D_1 \in { \bf \cal C}_1$. Meantime
$E {~and~} E_1 \in  { \bf \cal C}_2$. By Definition \ref{7.1.4},
$D \sim E$.

(3)  Obviously, $M_F \sim M_D$ for any $M_F \in [M_D]$. If $M_D
\not= 0$,
  then $M_F \not= 0$ by Lemma \ref {7.1.7} (1),. By Lemma \ref {7.1.7}(2), $ F \sim D$.
Thus  $M_F \in \{ M_E \mid E \in [D] \}$. If $M_D = 0$, by Lemma
\ref{7.1.7}(1), $M_F = 0 = M_D \in \{ M_E \mid E \in [D] \}$.

(4)  If $M_F \in [M_D]$, then $M_D \sim M_F$. If $D = 0$, then
$M_D = 0$. By Lemma \ref {7.1.7}(1), $ M_F = 0$. Thus $M_F = M_D =
0 \in$ \{$ M_E  \mid E \in [D]$ and $E$ is faithful to  $M$ \}. If
$ D \not= 0$, we have that
 $M_D \not= 0$ since  $D$ is faithful to  $M$.
By  Lemma \ref {7.1.7}(1), $M_F \not= 0$. By Lemma \ref{7.1.7}(3),
$M_F\in $ \{ $M_E  \mid E \in [D]$ and  $E$ is faithful to $M$ \}.
$\Box$

\begin{Theorem} \label{7.1.8}
Let $\{ { \bf \cal E}_{\alpha} \mid \alpha \in \bar{\Omega} \}$ be
all of the equivalence classes of   $C$. and $E_\alpha = \sum \{ E
\mid E \in { \bf \cal E}_\alpha \}$. Then

(1)  For any  $\alpha \in \bar{ \Omega }$, there exists a set $I_
\alpha$ and subclasses
 ${ \bf \cal E}(\alpha , i ) \subseteq { \bf \cal E}_{ \alpha }$
such that $\cup \{ { \bf \cal E}( \alpha,i) \mid i \in I_{ \alpha
} \} = { \bf \cal E}_{ \alpha }$ and $\{ M_{ { \bf \cal E}(
\alpha,i)}  \mid   \alpha \in \bar \Omega ,
i  \in I_ \alpha  \}$ \\
 is the set of the equivalence classes of  $M$
(they are distinct except for  zero equivalence class), where $M_{
{ \bf \cal E}( \alpha ,i)}$ denotes $ \{ M_D \mid  D \in { \bf
\cal E}( \alpha,i) \}$.

(2)  If  $M$ is a component faithfulness  $C$-comodule, then
 $ \{ M_{ { \bf \cal E}(\alpha , i)}  \mid
 \alpha \in \bar{ \Omega } , i \in I_{ \alpha } \}$
is the set  of the distinct equivalence  classes of  $M$.

(3) $$ M = \sum_{ \alpha \in \bar{ \Omega }}  \oplus M_{(E_{
\alpha })^{ (\infty)}} =  \sum_{ \alpha \in \bar{ \Omega }} \sum_{
i \in I_{ \alpha } }
 \oplus (M_{E( \alpha,i)})^{ (\infty)}
= \sum_{ \alpha \in \bar{ \Omega }}  \oplus (M_{E_{ \alpha }})^{
(\infty)}
$$
and for any
 $\alpha \in \bar{ \Omega }$,

$$ M_{(E_{ \alpha })^{ (\infty)}}
= \sum_{  i \in I_{ \alpha } }
 (M_{E( \alpha,i)})^{ (\infty)}
= ( M_{E_{ \alpha }})^{ (\infty)} $$ where $E( \alpha,i) = \sum \{
E  \mid  E  \in { \bf \cal E}( \alpha,i) \}$.

\end{Theorem}
 {\bf Proof}. (1) By Theorem \ref{7.1.3} and  Lemma \ref{7.1.7}(3) we can
immediately get part (1).

(2)  If  $M$ is a component faithfulness $C$-comodule, then ${ \bf
\cal C}(M)_1  =  { \bf \cal C}(M)_0 =  { \bf \cal C}_0  $. It
follows from Theorem \ref {7.1.3} and part (1)  that $ \{ M_{ {
\bf \cal E}(\alpha , i)}  \mid  \alpha \in \bar{ \Omega }, i \in
I_{ \alpha } \}$ consists of  all  the distinct equivalence
classes of $M$.

(3) We see that

 $M = M_C =
M_{ \sum \{ (E_{ \alpha})^{ (\infty)} \mid  \alpha \in \bar{
\Omega } \}}$
 (by  \cite{XF92})

$\supseteq \sum_{ \alpha \in \bar{ \Omega }} M_{{(E_{\alpha})}^{
(\infty)}}$ (by  Lemma \ref{7.1.6}(3))

$\supseteq \sum _{ \alpha \in \bar{ \Omega }} \sum _{ n = 0 }^{
\infty} M_{ \wedge ^{n+1}E_{ \alpha}}$ ( by Lemma \ref{7.1.6}(3))

$\supseteq \sum _{ \alpha \in \bar{ \Omega }} \sum _{ n = 0 }^{
\infty} {\wedge}^{n+1} M_{E_{ \alpha}}$ (by  Lemma \ref{7.1.6}(1))

$= \sum_{ \alpha \in \bar{ \Omega }} \sum _{ n = 0 } ^{ \infty}
 {\wedge} ^{n+1} M_{ \sum_{ i \in I_{ \alpha }}E( \alpha,i)}$  (by Theorem
\ref{7.1.8}(1) )

$\supseteq \sum _{ \alpha \in \bar{ \Omega }} \sum _{ n = 0 } ^{
\infty} {\wedge} ^{n+1} \sum _{ i \in I_{ \alpha }} M_{E( \alpha,
i)}$ (by Lemma \ref{7.1.6}(3))

$\supseteq \sum _{ \alpha \in \bar{ \Omega } }\sum _{ n = 0 } ^{
\infty}
 \sum_{ i \in I_{ \alpha }} {\wedge} ^{n+1} M_{E(\alpha,i)}$
( by  Lemma \ref{7.1.6}(5)(6))

$=  \sum _{ \alpha \in \bar{ \Omega }} \sum_{ i \in I_{ \alpha}}
 \sum _{ n = 0 } ^{ \infty} \wedge^{n+1}M_{E( \alpha,i)}$

$= \sum _{ \alpha \in \bar{ \Omega }} \sum_{ i \in I_{ \alpha }}
(M_{ E(\alpha,i)})^{ (\infty)}$

$=M$  (by \cite[(4.10) in Theorem4.15]{XSF94} and Lemma
\ref{7.1.6}(4)
 and part (1)).                                           \\
Thus
\begin {eqnarray}
 M = \sum_{ \alpha \in \bar{ \Omega }} M_{(E_{\alpha })^{ (\infty)}}
= \sum_{ \alpha \in \bar{ \Omega }} \sum_{ i \in I_{ \alpha }}
(M_{E( \alpha,i)})^{ (\infty)} \label {e.8.4}
\end {eqnarray}
and

$$M = \sum _{ \alpha \in \bar{ \Omega }} \sum_{ i \in I_{ \alpha }}\oplus
(M_{E( \alpha,i)})^{ (\infty)}$$ by \cite[Theorem 4.15]{XSF94} and
part (1). We see that
\begin {eqnarray*}
M_{E(\alpha,i)} \wedge M_{E(\alpha,i)}
&\subseteq & M_{(E_{\alpha})^{ (\infty)}} \wedge M_{(E_{\alpha})^{ (\infty)}} \\
&\subseteq & M_{(E_{\alpha})^{ (\infty)} \wedge (E_{\alpha})^{
(\infty)}}
\hbox { \ (by Lemma \ref{7.1.6}(1)) }  \\
&=& M_{(E_{\alpha})^{ (\infty)}} \hbox { (by \cite[ Proposition
2.1.1 ]  {HR74})},
\end {eqnarray*}
Thus $$M_{(E_{\alpha})^{ (\infty)}} \supseteq (M_{E( \alpha,i)})^{
(\infty)}$$  for any  $i  \in I_{\alpha }$ and   for any  $\alpha
\in \bar{ \Omega }$,  and
 \begin {eqnarray}
 M_{(E_{ \alpha})^{ (\infty)}} \supseteq
\sum_{ i \in I_{ \alpha }}(M_{E( \alpha,i)})^{ (\infty)} \label
{e.8.5}
\end {eqnarray}
If  $M_{(E_{\alpha})^{(\infty)}} \cap \sum_{\beta \in \bar{ \Omega
}, \beta \not= \alpha } M_{(E_{ \beta})^{ (\infty)}} \not= 0$,
then there exists a non-zero simple subcomodule $C^*x \subseteq
M_{(E_{ \alpha})^{ (\infty)}} \cap
 \sum_{ \beta \in \bar{ \Omega },
\beta \not= \alpha } M_{(E_{ \beta})^{ (\infty)}}$.
 By
\cite[Proposition 1.16]{XSF94}, $<x>$ is a minimal closed
subcomodule of $M$. By Lemma \ref{7.1.1}(3),
   $M_{(E_{ \alpha})^{ (\infty)}}$
is a closed subcomodule of $M$. By \cite[Lemma 3.3]{XSF94},
 there exists  $\gamma \in  \bar{ \Omega }$ with
$ \gamma \not= \alpha$ such that   $C^*x \subseteq M_{(E_{
\gamma})^{(\infty)}}$.
 Thus  $M_{(E_{ \gamma})^{(\infty)}} \cap M_{(E_{ \alpha })^{(\infty)}}
\not= 0$. By Lemma \ref {7.1.1} (7), $ (E_{ \gamma})^{(\infty)}
\cap (E_{ \alpha })^{(\infty)} \not= 0,$ which  contradicts
\cite{XF92}. Thus for any $\alpha \in \bar{ \Omega }$, $M_{(E_{
\alpha})^{ (\infty)}} \cap \sum_{ \beta \in \bar{ \Omega }, \beta
\not= \alpha } M_{(E_{ \alpha})^{ (\infty)}}= 0,$ which implies
that
$$M = \sum_{ \alpha \in \bar{ \Omega }} \oplus
 M_{(E_{ \alpha})^{ (\infty)}}.$$

It follows from equations ( \ref {e.8.4}) and (\ref {e.8.5} ) that
\begin {eqnarray}
                M_{(E_{ \alpha})^{ (\infty)}} &=&
\sum_{ i \in I_{ \alpha }}(M_{E( \alpha,i)})^{ (\infty)} \label
{e.8.6}
\end {eqnarray}
We see that
\begin {eqnarray*}
M_{E_{\alpha}} \wedge M_{E_{\alpha}} & \subseteq & M_{E_{\alpha}
\wedge E_{\alpha}}
 \hbox { \ (by Lemma \ref {7.1.6} (1)) } \\
 & \subseteq & M_{(E_{\alpha})^{ (\infty)}}.
 \end {eqnarray*}
 Thus
$$M_{(E_{ \alpha})^{(\infty)}} \supseteq (M_{E_{ \alpha}})^{(\infty)}
\supseteq \sum_{ i \in I_{ \alpha }}(M_{E( \alpha,i)})^{
(\infty)}$$ and   $$M_{(E_{ \alpha})^{(\infty)}}= (M_{E_{
\alpha}})^{(\infty)}$$ by relation ( \ref {e.8.6} ).
 This completes the proof. $\Box$

\begin {Definition} \label{7.1.9}  If  $M$
 can   be decomposed
into a direct sum  of two non-zero  weak-closed subcomodules, then
$M$ is called decomposable.      If $N$ is a subcomodule of $M$
and  $N$ contains exactly one non-zero minimal closed submodule,
then $N$  is said to be relative-irreducible.

\end {Definition}

\begin{Corollary} \label{7.1.10}   $C$ is a coalgebra.

(1)  If $C$ is   $\pi$-commutative, then every  $C$-comodule
 $M$ is also $\pi$-commutative;

(2)  If  $C$ can be decomposed into a direct sum of its
irreducible subcoalgebras, then every $C$-comodule
 $M$ can also be decomposed into a direct sum of its relative-irreducible
subcomodules;

(3)  If  $C$ is decomposable, then every component faithfulness
 $C$-comodule  $M$ is decomposable;

(4) If  $C$ is irreducible,  then every non-zero  $C$-comodule
$M$ is relative-irreducible;

(5) $C$ is irreducible iff every component faithfulness
$C$-comodule  $M$ is relative-irreducible.

\end{Corollary}
{\bf Proof}. (1) For any pair of non-zero closed subcomodules $N$
and $L$ of $M$, by Theorem \ref{7.1.3}, there exist  faithful
simple subcoalgebras $D$ and $E$ of $C$ to $M$  such that $N=M_D$
and $L =M_E$.
 By Lemma \ref{7.1.6}(2),
$ N \wedge L = M_D \wedge M_E = M_{D \wedge E} = M_{E \wedge D } =
M_E \wedge M_D = L \wedge N$. Thus $M$ is  $\pi$-commutative.

(2) Since  $C$ can be decomposed  into a direct sum of its
irreducible
 subcoalgebras,   every equivalence class of $C$
contains only one element by \cite{XF92}. By Theorem
\ref{7.1.8}(1), every equivalence class of $M$ also contains only
one element. Thus it follows from
 \cite [Theorem 4.18]{XSF94}  that
 $M$ can be decomposed into a direct sum of its relative-irreducible
 subcomodules.

(3)  If  $C$ is decomposable, then there are at least two non-zero
equivalence classes in  $C$. By   Theorem \ref{7.1.8}(2),
 there are at least two non-zero equivalence classes in $M$.
Thus  $M$ is decomposable.

(4)  If  $C$ is irreducible, then  there is only one non-zero
simple subcoalgebra of $C$ and there is at most one non-zero
minimal closed subcomodule in $M$ by Theorem \ref {7.1.3}.
Considering  $M \not= 0$, we have that $M$ is
relative-irreducible.

(5)  If  $C$ is irreducible, then every component faithfulness
$C$-comodule  $M$ is relative-irreducible
 by  Corollary \ref{7.1.10}(4).
Conversely, let $M = C$  be the regular  $C$-comodule. If  $D$ is
a non-zero simple subcoalgebra of  $C$, then
 $ 0 \not= D \subseteq M_D$. Thus $M$ is a component faithfulness
$C$-comodule. Since  $M$
 is a relative-irreducible $C$-comodule by assumption,
  there is only one non-zero
minimal closed subcomodule in  $M$ and so
 there is also only  one non-zero simple subcoalgebra in $C$
 by Theorem \ref {7.1.3}.
Thus  $C$ is irreducible. $\Box$

\begin{Lemma} \label{7.1.11} Let $N$ be a
$C$-subcomodule of $M$ and $\emptyset \not= L \subseteq M$. Then

(1)  $C^* \cdot L = C(M)^* \cdot L$; {~}{~} $N^{ \bot C^* \bot M }
= N^{ \bot C(M)^* \bot M }$;

(2)  $N$ is a (weak-) closed  $C$-subcomodule iff  $N$ is a
(weak-)closed $C(M)$-subcomodule;

(3)  $N$ is a minimal closed  $C$-subcomodule iff $N$ is a minimal
closed $C(M)$-subcomodule;

(4) $N$ is an indecomposable  $C$-subcomodule iff $N$ is an
indecomposable  $C(M)$-subcomodule;

(5)  $N$  is a relative-irreducible  $C$-subcomodule iff $N$ is a
relative-irreducible  $C(M)$-subcomodule.
\end{Lemma}
{\bf Proof}. (1) Let $C=C(M) \oplus V,$ where $V$ is a subspace of
$C$. If  $f \in V^*$, then  $ f \cdot L = 0$. Thus
  $C^* \cdot L  = (C(M)^* + V^* ) \cdot L = C(M)^* \cdot L$.
We now show the second equation. Obviously,
 $$N^{ \bot C^* \bot M } \subseteq
N^{ \bot C(M)^* \bot M }.$$ Conversely,  for any $x \in N^{ \bot
C(M)^* \bot M }$ and $f \in N^{ \bot C^* }$, there exist  $f_1 \in
C(M)^*$ and  $f_2 \in V^*$ such that
 $f = f_1 + f_2$.     Obviously,
  $f \cdot x = f_1 \cdot x = 0.$  Thus
$x \in N^{ \bot C^* \bot M }$ and $N^{ \bot C(M)^* \bot M }
\subseteq N^{ \bot C^* \bot M }$.
 Therefore
 $$N^{ \bot C^* \bot M } = N^{ \bot C(M)^* \bot M }.$$

(2)  If  $N$ is a weak-closed  $C$-subcomodule, then $(C^* \cdot
x)^{ \bot C^* \bot M } \subseteq N $ for any $x \in N$. We see
that
\begin {eqnarray*}  (C^* \cdot x)^{ \bot C^* \bot M }
&=& (C(M)^* \cdot x)^{ \bot C^* \bot M } \hbox { ( \ by  Lemma \ref{7.1.11}(1))  } \\
&=& (C(M)^* \cdot x)^{ \bot C(M)^* \bot M } \hbox { \ ( \ by Lemma
\ref{7.1.11}(1)) }
\end {eqnarray*}
Thus
 $(C(M)^* \cdot x)^{ \bot C(M)^* \bot M } \subseteq N$ and so
 $N$ is a weak-closed
 $C(M)$-subcomodule.
Conversely, if  $N$ is a weak-closed  $C(M)$-subcomodule.
Similarly, we can show that  $N$ is a weak-closed
$C$-subcomodule. We now show the second assertion.  If $N$ is a
closed $C$-subcomodule, then  $N^{ \bot C^* \bot M } =  N$ and
 $N= N^{ \bot C^* \bot M } = N^{ \bot C(M)^* \bot M } $ by part (1),
 which implies that
 $N$ is a closed  $C(M)$-subcomodule. Conversely,
if  $N$ is a closed  $C(M)$-subcomodule,
 similarly, we can  show that $N$ is a closed $C$-subcomodule.

Similarly the others can  be proved. $\Box$

\begin{Proposition} \label{7.1.12}  Let
every simple subcoalgebra in $C(M)$ be faithful to  $M$.

(1)  If $M$ is an indecomposable $C$-comodule, then $C(M)$ is also
an indecomposable subcoalgebra.

(2)  $M$ is a  relative-irreducible  $C$-comodule iff $C(M)$ is an
irreducible subcoalgebra.
\end{Proposition}
 {\bf Proof}.
(1)  If  $M$ is an indecomposable  $C$-comodule, then $M$ is an
indecomposable  $C(M)$-comodule by Lemma \ref{7.1.11} and $M$ is a
component faithfulness $C(M)$-comodule. By Corollary
\ref{7.1.10}(3),  $C(M)$ is indecomposable.

(2)  If  $M$ is a relative-irreducible $C$-comodule, then
  $M$ is a relative-irreducible  $C(M)$-comodule
by Lemma \ref{7.1.11}(5)  and
  $C(M)$ is irreducible. Conversely,
if  $C(M)$ is irreducible, then $M$ is a relative-irreducible
$C(M)$-comodule  by Corollary \ref{7.1.10}(4) and so $M$ is a
relative-irreducible $C$-comodule by Lemma \ref{7.1.11}(5). $\Box$

\begin{Definition} \label{7.1.13}
If
 $D \sim E$ implies  $M_D \sim M_E$ for any simple subcoalgebras
$D$ and $E$ in $C(M)$, then $M$ is called a $W$-relational
hereditary
 $C$-comodule.
\end {Definition}

 If $M$ is a  $W$-relational hereditary $C$-comodule, then
 ${ \cal C}(M)_1 = {\cal C}(M)_0$.
In fact, by Lemma \ref {7.1.1}(5) , ${ \cal C}(M)_1 \subseteq
{\cal C}(M)_0$. If ${ \cal C}(M)_1 \not= {\cal C}(M)_0$, then
there exists $ 0 \not= D \in { \cal C}(M)_0$ such that  $M_D = 0$.
Since  $M_D \sim M_0$ and  $M$ is  $W$-relational hereditary,  we
have that $D \sim 0$ and $D = 0$ by Lemma \ref {7.1.7}(1). We get
a contradiction. Thus
 ${ \cal C}(M)_1 = {\cal C}(M)_0$.

Obviously, every $\pi$-commutative comodule is  $W$-relational
hereditary. If $C$ is $\pi$-commutative, then $M$ is
$\pi$-commutative by Corollary \ref {7.1.10}(1) and  every
$C$-comodule $M$ is $W$-relational hereditary.
 Furthermore, $M$ is also a component faithfulness $C(M)$-comodule.

\begin{Proposition} \label{7.1.14}
Let  the notation be the same as in Theorem \ref{7.1.8}. Then the
following conditions are equivalent.

(1)  $M$ is  $W$-relational hereditary.

(2)  For any  $\alpha \in  \bar { \Omega }$,
 there is at most one non-zero equivalence class in
 $\{ M_{ { \bf \cal E}( \alpha, i)} \mid i \in I_{ \alpha } \},$
and   ${ \cal C}(M)_1 = {\cal C}(M)_0.$

(3)    For any   $ D \in { \bf \cal C}(M)_0$,

 $[M_D] = $\{ $M_F \mid F \in [D]$  and  $ F \subseteq C(M)$ \}

(4)  For any  $D$ and  $E$ $\in { \bf \cal C}(M)_0$ ,
 $M_D \sim M_E$ iff $D \sim E$.

(5)  For any   $\alpha \in  \bar { \Omega }$, $M_{(E_{\alpha})^{
(\infty) }}$ is indecomposable, and  ${ \cal C}(M)_1 = {\cal
C}(M)_0.$

(6) For any  $ \alpha \in  \bar { \Omega } $, $(M_{E_{\alpha}})^{
(\infty) }$
 is indecomposable, and
  ${ \cal C}(M)_1 = {\cal C}(M)_0;$

(7) $ M = \sum_{ \alpha \in \bar{ \Omega }}  \oplus (M_{E_{ \alpha
}})^{ (\infty)}$ is its weak-closed indecomposable decomposition,
and
 ${ \cal C}(M)_1 = {\cal C}(M)_0;$

(8)  $ M = \sum_{ \alpha \in \bar{ \Omega }}  \oplus M_{(E_{
\alpha })^{ (\infty)}}$ is its weak-closed indecomposable
decomposition, and
 ${ \cal C}(M)_1 = {\cal C}(M)_0.$
 \end{Proposition}
{\bf Proof}. We prove it along with the following lines:
 $(1) \Longrightarrow  (4)
 \Longrightarrow  (3) \Longrightarrow  (1)$
$(3) \Longleftrightarrow  (2) \Longleftrightarrow  (5)
 \Longleftrightarrow  (6)  \Longleftrightarrow  (7)  \Longleftrightarrow  (8)$.

$(1) \Longrightarrow  (4)$ It follows from the discussion above
that
$${ \cal C}(M)_1 = {\cal C}(M)_0.$$
If  $D$ and  $E \in {  \cal C}(M)_0$
 and $M_D \sim M_E $, then   $ D \sim E  $ by  Lemma \ref {7.1.7} (2). If
$D \sim E$ and  $D$ and $E \in {  \cal C}(M)_0$, then $M_ D \sim
M_E$.

$(4) \Longrightarrow (3)$  Considering  Lemma \ref{7.1.7}(4) and
 ${ \cal C}(M)_1 = {\cal C}(M)_0$,
 we only need to show that
$$ \{ M_F \mid F \in [D] \hbox { \  and \ } F \subseteq C(M) \}
 \subseteq [M_D].$$ For any $F \in [D]$ with $F \subseteq C(M)$,
$  F \sim D$ and  $M_F \sim M_D$ by part (4), which implies that
$M_F  \in  [M_D]$.

$(3) \Longrightarrow (1)$  It is trivial.

$(5) \Longleftrightarrow (6)
 \Longleftrightarrow (7)
 \Longleftrightarrow (8)$
 It follows from  Theorem \ref{7.1.8}(3).

$(5) \Longleftrightarrow (2)$ By Theorem \ref{7.1.8} (3),
$$M_{(E_{\alpha})^{ (\infty )}}= \sum_{ i \in I_{ \alpha}}
(M_{E(\alpha,i)})^{ (\infty) }.$$ Thus part (2) and part (5) are
equivalent.

$(2) \Longrightarrow  (3)$ It follows from  Lemma \ref {7.1.7}
(4).

 $(3) \Longrightarrow  (2)$ If there are two non-zero
equivalence classes  $M_{ { \bf \cal E}( \alpha, i_1)}  \not= 0 $
and $ M_{ { \bf \cal E}( \alpha, i_2)} \not= 0$ in $\{ M_{ { \bf
\cal E}( \alpha, i) } \mid i \in I_{ \alpha } \}$, then there
exist
 $D_1 \in   { \bf \cal E}( \alpha, i_1)$ and
$D_2 \in  { \bf \cal E}( \alpha, i_2)$ such that $M_{D_1} \not= 0$
and $M_{D_2} \not=0$. Let  $D = D_1$. Since  $M_{D_1}$ and
$M_{D_2} \in$
 \{$M_F \mid F \in [D]$  and $F \subseteq C(M)$ \} and
 $M_{D_2} \not\in [M_D]
 = M_{ { \bf \cal E}( \alpha, i_1)}$, this contradicts part (3).
Thus there is at most one non-zero equivalence class in
 $\{ M_{ { \bf \cal E}( \alpha, i)} \mid i \in I_{ \alpha } \}$.
$\Box$

\begin{Proposition} \label{7.1.15} If
  $M$ is a full,   $W$-relational hereditary $C$-comodule,
  then

(1) $M$ is indecomposable iff $C$ is indecomposable;

(2) $M$ is relative-irreducible iff $C$ is irreducible.

(3) $M$ can be decomposed into a direct sum of
 its weak-closed  relative-irreducible subcomodules iff
$C$ can be decomposed into a direct sum of its irreducible
subcoalgebras.

(4) $M$ is  $\pi$-commutative iff $C$ is  $\pi$-commutative.

\end{Proposition}
 {\bf Proof}. Since  $M$ is a full $C$-comodule,  $C(M) =C$.
Since  $M$ is  $W$-relational hereditary,  ${ \cal C}(M)_1 = {\cal
C}(M)_0 = { \cal C}_0$. Thus  $M$ is a component faithfulness
$C$-comodule.

(1) If  $M$ is indecomposable, then $C$ is indecomposable by
Proposition \ref {7.1.12}(1). Conversely, if  $C$ is
indecomposable, then there is at most one non-zero equivalence
class in  $C$. By Proposition \ref {7.1.14}(2), there is at most
one non-zero equivalence class in  $M$. Thus
 $M$ is indecomposable.

(4) If  $M$  is  $\pi$-commutative, then there is only one element
in every  equivalence class of $M$. By Proposition \ref
{7.1.14}(4) and Theorem \ref{7.1.3}, there is only one element in
every  equivalence class of   $C$.
 Thus  $C$ is  $\pi$-commutative by \cite {XF92}.
Conversely, if $C$  is $\pi$-commutative, then $M$ is
$\pi$-commutative by Corollary \ref {7.1.10}(1).

 (2) It follows from the above discussion
 and Proposition \ref {7.1.12}.

$(3) \Longleftrightarrow (4)$
 By  \cite[Theorem3.8 and Theorem 4.18]{XSF94} and \cite{XF92},
 it is easy to check that (3) and (4) are equivalent. $\Box$

\begin {Proposition} \label {7.1.16}
 If  $M_D
 \wedge M_E =M_E \wedge M_D$ implies  $D \wedge E = E \wedge D$
for any simple coalgebras $D$ and $E$ in $C(M)$, then  $M$ is
$W$-relational hereditary.
\end {Proposition}
{\bf Proof}.  Let $ D \sim E$. For any  pair of subclasses ${ \bf
\cal M}_1$ and  ${ \bf \cal M}_2$ of ${ \bf \cal M}_0$ with
 ${ \bf \cal M}_1 \cap { \bf \cal M}_2
= \emptyset $ and   ${ \bf \cal M}_1  \cup { \bf \cal M}_2 = { \bf
\cal M}_0$, if  $ M_D \in { \bf \cal M}_1$ and $M_E \in { \bf \cal
M}_2$, let
 ${ \bf \cal C}_1  =
\{ F \in { \bf \cal C}_0  \mid M_F \in { \bf \cal M}_1  \}$ and
 ${ \bf \cal C}_2  = \{ F \in { \bf \cal C}_0  \mid M_F \in { \bf \cal M}_2  \}$.
By  Theorem \ref{7.1.3},
 $ { \bf \cal C}_0 = { \bf \cal C}_1 \cup { \bf \cal C}_2$ and
 ${ \bf \cal C}_1 \cap { \bf \cal C}_2 =
\emptyset $. Obviously, $D \in  { \bf \cal C}_1 $ and $ E \in {
\bf \cal C}_2$. By Definition \ref {7.1.4}, there exist $D_1 \in {
\bf \cal C}_1$ and $E_1 \in { \bf \cal C}_2 $ such that
 $D_1 \wedge  E_1 \not=   E_1 \wedge  D_1$. By the assumption condition,
we have that
 $M_ {D_1} \wedge M_{E_1} \not=  M_{E_1} \wedge M_{D_1}$. Thus
$M_D \sim M_E$, i.e. $M$ is  $W$-relational hereditary. $\Box$

\begin{Proposition} \label{7.1.17}
Let  $M=C$ as a right $C$-comodule. Let $N=D$ and $L=E$ with $N$
and $L$ as subcomodules of $M$ with $D$ and $E$ as right coideals
of $C$. Let $X$ be an ideal of $C^*$ and $F$ be a subcoalgebra of
$C$. Then:

(1) $ X^{ \bot C } = X^{ \bot M}$;{~}{~}
  $ N^{ \bot C^* \bot M}  = C(N)$;     {~}{~}
$C(M_F) = F$;

(2)  $N$ is a closed subcomodule of  $M$ iff $D$ is a subcoalgebra
of $C$;

(3)   $N$ is a closed subcomodule iff  $N$ is a weak-closed
subcomodule of $M$;

(4) $N$ is a minimal closed subcomodule of $M$ iff  $D$  is  a
simple subcoalgebra of $C$;

(5) When  $N$ and  $L$ are closed subcomodules, $N \wedge _M L  =
D \wedge _C E$, where  $ \wedge _M$ and  $\wedge _C $ denote
wedge  in comodule $M$ and  in coalgebra $C$ respectively;

(6)   $M$ is a full and  $W$-relational hereditary $C$-comodule
and a component faithfulness $C$-comodule.

(7) The weak-closed indecomposable decomposition of $M$
 as a $C$-comodule and
the indecomposable decomposition of $C$ as coalgebra are the same.

\end{Proposition}
{\bf Proof}. (1)  If $ x \in  X^{ \bot M}$, then $  f \cdot x = 0$
for any  $ f \in X$.
 Let  $\rho (x) = \sum_{ i = 1 } ^{ n} x_i \otimes c_i $ and
$ x_1$, $\cdots $, $x_n$ be linearly independent. Thus  $ f(c_i) =
0$ and  $i = 1$, $ \cdots$ ,$ n$.
 Since
$x = \sum _{ i = 1 }^{n} \epsilon (x_i)c_i$,
 $f(x) = 0$, which implies  $x \in X^{ \bot C}$.

Conversely, if  $x \in X^{ \bot C}$, we have that $\rho (x) \in
  X^{ \bot C } \otimes  X^{ \bot C}$
since  $X^ {\bot C}$  is subcoalgebra of $C$. Thus  $f \cdot x =
0$ for any  $f \in X$, which implies $ x \in X^{ \bot M}$. Thus $
X^{ \bot M}  =  X^{ \bot C }$. By Lemma \ref{7.1.1}(1),  $ N^{
\bot C^* \bot M}  = C(N)^{ \bot C^* \bot M}$. By part (1),  $
C(N)^{ \bot C^* \bot M}
 = C(N)^{ \bot C^* \bot C}$.
Thus $ N^{ \bot C^* \bot M}  =C(N)^{ \bot C^* \bot C}  = C(N)$.

Finally, we show that  $C(M_F) = F$. Obviously,  $ C(M_F)
\subseteq F$. If we view $F$ as a $C$-subcomodule of $M$, then
  $F \subseteq C(M_F)$. Thus  $F = C(M_F)$.

(2)  If  $N$ is a closed subcomodule of  $M$, then $  N^{ \bot C^*
\bot M}  = N$. By Proposition \ref {7.1.17}(1),
  $ N^{ \bot C^* \bot M}  = C(N) = N$. Thus
 $\rho (N) = \Delta (D)\subseteq
N \otimes N = D \otimes D$, which implies  that $D$ is a
subcoalgebra of $C$. Conversely, if  $D$ is a subcoalgebra of $C$,
then
 $\rho (N) \subseteq N \otimes N$.
Thus  $C(N) \subseteq D = N$. By Proposition \ref{7.1.17}(1),
   $N^{ \bot C^* \bot M}  = C(N) \subseteq N$. Thus
   $N^{ \bot C^* \bot M}  = N$, i.e. $N$ is closed.

(3) If  $N$ is a closed subcomodule, then  $N$ is weak-closed.
Conversely, if  $N$ is weak-closed, then $ <x> \subseteq N$
 for  $x \in N$. Since  $<x>$
is a closed subcomodule of  $M$,
 $ <x>$ is
subcoalgebra of $C$ if let $<x>$ with structure of coalgebra $C$.
This shows that $ \Delta (x) \in <x> \otimes <x> \subseteq D
\otimes D$. Thus  $D$ is a subcoalgebra of $C$. By Proposition
\ref {7.1.17}(2), $N$ is a closed subcomodule of $M$.

(4) It follows from  part (2).

(5) We only need to show that
$$  ( N^{ \bot C^* } L^{ \bot C^*})^{ \bot M }=
(  D^{ \bot C^* }E^{ \bot C^*})^{\bot C}.$$ Since  $N$ and  $L$
are closed subcomodules,
  $N^{ \bot C^* } = C(N)^{ \bot C^*} = D^{ \bot C^* }$ and
$L ^{ \bot C^* }  = C(L)^{ \bot C^* } = E^{\bot C^*}$ by
Proposition \ref {7.1.17}(1) and Lemma \ref{7.1.1}(1).
 Thus we only need to show that

$(D^{ \bot C^*} E^{ \bot C^*})^{ \bot M}  =
(D^{ \bot C^* } E^{ \bot C^*})^{ \bot C }$.  \\
The above formula follows from Proposition \ref {7.1.17}(1).

(6)  By the proof of Corollary \ref {7.1.10}(5), we know that $M$
is  a component faithfulness $C$-comodule. Let $\{m_{ \lambda }
\mid  \lambda \in \Lambda \}$ be a basis of $M$. For any   $c \in
C$, $\Delta(c) =
        \sum m_{ \lambda } \otimes c_{\lambda }$,
by the definition of  $C(M)$, $ c_{\lambda } \in C(M)$. Since
$c= \sum \epsilon(m_{ \lambda }) c_{\lambda } \in C(M)$, $ C
\subseteq C(M)$, i.e. $ C = C(M)$. Consequently, it follows from
part (4)(5) that  $M$ is  $W$-relational hereditary.

(7)  Since  $M$ is  $W$-relational hereditary,
 $ M = \sum_{ \alpha \in \bar{ \Omega }}
 \oplus M_{(E_{ \alpha })^{ (\infty)}}$
is a weak-closed indecomposable decomposition  of $M$ by
Proposition \ref{7.1.14}. By part (1) (3), $ C(M_{(E_{ \alpha })^{
(\infty)}}) =(E_{ \alpha })^{ (\infty)}  =
 (M_{(E_{ \alpha })^{ (\infty)}})^{ \bot C^* \bot M}
=  M_{(E_{ \alpha })^{ (\infty)}}$. Thus   $ M = \sum_{ \alpha \in
\bar{ \Omega }} \oplus M_{(E_{ \alpha })^{ (\infty)}} = \sum_{
\alpha \in \bar{ \Omega }}  \oplus (E_{ \alpha })^{ (\infty)}$. By
\cite{XF92},
 $ \sum_{ \alpha \in \bar{ \Omega }}  \oplus (E_{ \alpha })^{ (\infty)}
= C$ is a indecomposable decomposition  of $C$. Thus the
weak-closed indecomposable decomposition  of $M$
 as a $C$-comodule and
the indecomposable decomposition  of $C$ as coalgebra are the
same.
        This completes the proof. $\Box$

By Lemma \ref{7.1.6}(4) and Proposition \ref {7.1.17}, $C$ is
cosemisimple iff every $C$-comodule $M$ is cosemisimple.

\section { The coradicals of comodules}
\begin{Proposition} \label{7.2.1}
   Let  $M$ be a  $C$-comodule,  $J$ denote the
Jacobson radical of $C^*$ and   $r_j(C(M)^*)$ denote
 the Jacobson radical of $C(M)^*.$

(1)   $$M_0 = (r_j(C(M)^*))^{ \bot M } = Soc_{C^*}M \hbox { \ \
and \ \ } r_j(C(M)^*)= M_0^{ \bot C^*};$$

(2) If we view    $C$  as  a right  $C$-comodule, then
$$C_0 = Soc_{C^*}C = J^{ \bot C }
= \sum \{ D \mid D \hbox { \ is a minimal right coideal of \ }  C
\}.$$

\end{Proposition}
 {\bf Proof}.  (1)  We first show that
$$ M_0^{ \bot C^*}
\subseteq J$$ when  $M$ is a full  $C$-comodule. We only need to
show that  $ \epsilon - f$  is invertible in
 $C^*$ for any  $f \in (M_0)^{\bot C^*} $.
Let  $I = (M_0)^{ \bot C^*}$. For any natural number $n$, $f^{n+1}
\cdot (M_0)^{(n)} =0$ since
 $f^{n+1} \in I^{n+1}$,
 where $(M_0)^{(n)} = \wedge ^{n+1}M_0$.
 Thus
\begin {eqnarray}
f^{n+1}(C((M_0)^{(n)})) &=& 0 \label {e.18.1}.
\end {eqnarray}
 by  Lemma \ref {7.1.1}(1) and

\begin {eqnarray}      M=(M_0)^{ (\infty)}
&=& \cup \{ (M_0)^{(n)} \mid n = 0, 1, \cdots \}
 \label {e.18.2}
 \end {eqnarray} by \cite[Theorem 4.7]{XSF94}.
We now show that
\begin {eqnarray}
C  &=& \cup \{ C((M_0)^{(n)}) \mid n = 0, 1, \cdots \} \label
{e.18.3}
\end {eqnarray}
Since  $(M_0)^{(n)} \subseteq (M_0)^{(n+1)},$  we have that there
exists a basis  $ \{ m_\lambda \mid  \lambda \in \Lambda  \}$ of
$M$ such that for every given natural number  $n$ there exists a
subset of $ \{ m_\lambda \mid  \lambda \in \Lambda  \},$ which is
a basis of $(M_0)^{(n)}$. For any $c \in W(M)$, there exists  $m
\in M$ with $\rho(m) = \sum m_{ \lambda } \otimes c_{\lambda }$
such that
   $c_{\lambda _0 } = c $ for some $\lambda_0 \in \Lambda $.
  By equation (\ref {e.18.2}), there exists a natural number $n$ such that
   $m \in (M_0)^{(n)},$ which implies that
$c \in C((M_0)^{(n)})$ and $C(M) \subseteq  \cup _{0}^{\infty}
C((M_0)^{(n)})$. Thus
$$C= \cup \{ C((M_0)^{(n)}) \mid n = 0, 1, \cdots \}.$$
Let
$$g_n = \epsilon + f + \cdots + f^n , n = 1, 2, \cdots.$$

For any $c \in C$, there exists a natural number $n$
 such that
$c \in  C((M_0)^{(n)})$ by relation (\ref {e.18.3}). We define
\begin {eqnarray}
g(c) = g_n(c)\ . \label {e.18.4}
\end {eqnarray}
Considering relation (\ref{e.18.1}), we have  that $g$ is
well-defined. Thus  $g \in C^*$.

We next show that  $g$ is an inverse of  $ \epsilon -f $ in $C^*$.
For any $c \in C^*$,  there exists a natural number  $n$ such that
$c \in C((M_0)^{(n)})$ by relation (\ref {e.18.3}). Thus  $\Delta
(c) \in C((M_0)^{(n)}) \otimes C((M_0)^{(n)})$. We see that
\begin {eqnarray*}
g*( \epsilon -f)(c) &=& \sum g(c_1)(\epsilon -f)(c_2) \\
 &=& \sum g_n(c_1)(\epsilon -f)(c_2) \hbox { (by relation (\ref {e.18.4} )) } \\
 &=& (g_n *(\epsilon -f))(c)   \\
 &=& (\epsilon -f^{n+1})(c)  \\
 &=& \epsilon (c) \hbox { ( by relation  (\ref {e.18.1})) }.
 \end {eqnarray*}
Thus  $g*(\epsilon -f) = \epsilon $. Similarly, $(\epsilon -f)*g =
\epsilon $. Thus  $\epsilon -f$ has an inverse in  $C^*$.

        We next show that  $$ M_0^{ \bot C^*} = J$$ when $C(M)=C$.
It follows from Lemma \ref {7.1.1}(1) that

\begin {eqnarray*} (M_0)^{ \bot C^* }
&=& \cap \{ N^{\bot C^*} \mid N \hbox { \  is
a minimal closed subcomodule of \ } M \} \\
&=& \cap \{ C(N)^{ \bot C^*} \mid N \hbox { is a minimal closed
subcomodule
 of } M \}       \\
& \supseteq & J  \hbox { (  \ by \cite [Proposition 2.1.4] {HR74}
and
  Proposition \ref {7.1.2}) }
\end {eqnarray*}
Thus $M_0^{ \bot C^*} = J$

We now show that
 $$r_j(C(M)^*) = M_0 ^ {\bot C^*} \hbox { \ \ \ and \ \ \ }
 (r_j(C(M)^*))^{\bot M} = M_0$$
 for any $C$-comodule $M$.
If $M$ is  a $C$-comodule,  then $M$ is full $C(M)$-comodule and
so   $$(r_j(C(M)^*) = M_0^{\bot C^*}.$$ By  \cite [Proposition
4.6] {XSF94}, $M_0$ is closed. Thus   $((r_j(C(M)^*))^{\bot M} =
M_0$.

 Finally, we  show that   $M_0 = Soc_{C^*}M$ for any $C$-comodule $M$. By \cite
[Proposition 1.16]{XSF94}, $Soc_{C^*}M \subseteq M_0$. Conversely,
if  $x \in M_0$, let   $N = C^*x$. By Lemma \ref{7.1.6}(4),
 $M_{C_0} = M_0$. Thus
 $N$ is a finite dimensional  $C_0$-comodule. By  \cite [Theorem14.0.1]{Sw69a},
$N = N_1 \oplus \cdots \oplus N_n$ and
 $N_i$ is a simple  $C_0^*$-submodule. By  \cite [Lemma 1.1]{XSF94},
$N_i$ is a  $C_0$-subcomodule. Thus  $N_i$ is a $C$-comodule. By
\cite[Lemma 1.1]{XSF94}, $N_i$ is a  $C^*$-submodule. If  $L$ is
a non-zero $C^*$-submodule of $M$ and $L \subseteq N_i$, then  $L$
is also a  $(C_0)^*$-submodule. Thus
 $L = N_i$. This shows     that
 $N_i$ is also a simple  $C^*$-submodule. Thus  $N_i
\subseteq Soc_{C^*}M$ and  $N \subseteq Soc_{C^*}M$.
 It follows from the above proof that $M_0 = Soc_{C^*}M$.

(2) It follows from  Proposition \ref{7.1.17} and part (1).

\end {document}